\documentclass[10pt]{article}

\usepackage{amsmath,amssymb}

\newcommand{\Z}{\mathbb{Z}}

\newtheorem{theorem}{Theorem}
\newtheorem{lemma}[theorem]{Lemma}
\newtheorem{proposition}[theorem]{Proposition}

\DeclareMathOperator{\var}{Var}
\begin{document}
\title{A Random Walk with Collapsing Bonds and Its Scaling Limit}
\author{Majid Hosseini \\ and Krishnamurthi Ravishankar \\
Department of Mathematics,\\
 State University of New York at New Paltz\\
 1 Hawk Drive. Suite 9,\\
  New Paltz, NY 12561-2443\\
  {\ttfamily hosseinm@newpaltz.edu}\\
     {\ttfamily ravishak@newpaltz.edu}}
    \maketitle
\begin{abstract}
We introduce a new self-interacting random walk on the integers in
a dynamic random environment and show that it converges to a pure
diffusion in the scaling limit. We also find a lower bound on the
diffusion coefficient in some special cases. With minor changes
the same argument can be used to prove the scaling limit of the
corresponding walk in $\Z^d$.
\end{abstract}

In this note we introduce a self-interacting random walk in a
dynamic random environment, prove that it is recurrent and find
its scaling limit. The environment evolves in time in conjunction
with the random walk and is non-markovian;  however we will show
that the walk only remembers the recent past.
 Various models of self-interacting random walks and random walks in dynamic random environments
 have been studied recently (see for instance \cite{BandZeit06,
BenjaminiWilson03,BoldrighiniMinlosPellegrinotti97,BoldrighiniMinlosPellegrinotti2000,Davis99,
DavisVolkov02}  and references therein).

 Consider a particle performing a continuous-time nearest
neighbor symmetric random walk on the integers lattice. We assume
that the particle is initially at the origin. The times between
successive jumps are independent exponential random variables with
rate $\lambda$. Anytime the particle jumps over the bond
connecting two neighboring lattice sites, there is a probability
$0< p\leq  1$ that the bond connecting the sites breaks. The
particle is not able to jump over that bond until that bond is
repaired. If at the time the particle attempts to jump, one of the
bonds neighboring the particle is broken, the particle jumps over
the  other bond with probability one. If both bonds neighboring
the particle are broken, the particle can't jump when it attempts
to do so. The repair times of bonds are independent exponential
random variables with rate $\mu\in (0,\infty)$. Initially, there
are no broken bonds.

Let $X(t)$ denote the position of the particle at time $t\in
[0,\infty)$. We assume that the jump times are such that almost
surely, $X(t)$ is continuous from right and has a left limit at
all $t$. We will show that for all values of $p$, $X(t)$ is a
recurrent process and its scaling limit is a pure diffusion. We
will find a lower bound for the diffusion coefficient in some
special cases.

\begin{theorem}\label{th:recurrence}
The process $\{X(t)\}_{t\geq 0}$ is recurrent.
\end{theorem}
Theorem~\ref{th:recurrence} is a consequence of the following.
\begin{lemma}\label{lem:healing bonds} With probability 1, for any
$s\geq 0$, there exists $t\geq s$ such that at time $t$, there are
no broken bonds.
\end{lemma}
{\it Proof } The bonds break at a rate $\lambda p$, until the
particle is trapped at a lattice site. At that time no more bonds
break, until the particle is free to move again. Each bond is
repaired at a rate $\mu$, independent of the other bonds. Thus,
$b_t$,  the number of broken bonds at time $t$, is less than or
equal to $Q_t$, where $\{Q_t\}_{t\geq 0}$ is an $M/M/\infty$ queue
with incoming traffic rate $\lambda p$ and service rate $\mu$.

  It is straightforward to show that this
queue is a recurrent Markov chain. Furthermore if, starting from
zero customers, $T$ is the recurrence time of this queue back to
zero customers,  then $T$ is the sum of the idle time (which is an
exponential random variable with rate $\lambda p$) and the busy
period of this queue, and these two times are independent of each
other. The busy period has a Laplace transform (see
\cite{Shanbhag}).
 Thus $T$ is finite almost surly, has finite moments of all orders, and
$E(T)=\exp(\lambda p/\mu)/(\lambda p)$.  Note that since $b_t\leq
Q_t$, when $Q_t=0$ we also have $b_t=0$. Therefore, assuming there
are no broken bonds at the staring time, if $\tau$ is the
recurrence time of $b_t$ back to zero, then $\tau\leq T$. Thus
$\tau$ is finite almost surly, has finite moments of all orders,
and $E(\tau)\leq \exp(\lambda p/\mu)/(\lambda p)$. Hence, $0$ is a
recurrent state for $b_t$. Therefore, for any $s\geq 0$, there is
a $t\geq s$ such that there are no broken bonds at time $t$. This
completes the proof of the Lemma.

Put $\tau_0=0$ and if $\tau_i$ has been defined, put
\begin{equation*}\sigma_{i+1}=\inf\{t>\tau_i:b_t>0\}
\end{equation*} and
\begin{equation}\label{eq:taus}
\tau_{i+1}=\inf\{t>\sigma_{i+1}:b_t=0\}.
\end{equation}
The following is immediate from the proof of
Lemma~\ref{lem:healing bonds}.
\begin{lemma}\label{lem:taus}
The random variables $\tau_{i+1}-\tau_i$, $i\geq 0$ are i.i.d,
finite almost surely,  and have finite moments of all orders.
\end{lemma}

Let $N(t)$ be a Poisson process with parameter $\lambda$ which
increases by 1 anytime the particle makes an attempt to jump. Note
that there are times that the particle is trapped at a site and an
attempt to jump fails. We increase $N(t)$ by 1 in such instances
as well.

 {\it Proof of Theorem 1}

 By Lemma~\ref{lem:taus}, $\tau_1$ is finite almost surely
 and has finite moments of all orders. Since there are no broken bonds at times $\tau_i$, $i\geq 0$, the
sequence $\{X_{\tau_{i+1}}-X_{\tau_i}\}_{i\geq 0}$ is an i.i.d
sequence. Note that for any $t\geq 0$, $|X(t)|\leq N(t)$. Thus,
$E|X_{\tau_1}|\leq \lambda  E(\tau_1) \leq \exp(\lambda p/\mu)/p$.
By symmetry considerations, $E(X_{\tau_1})=0$. Thus,
$X_{\tau_n}=\sum_{k=0}^{n-1}(X_{\tau_{k+1}}-X_{\tau_k})$ is a sum
of mean zero, i.i.d. random variables. By standard arguments for
recurrence of sums, $P(X_{\tau_n}=0 \text{ i.o. } n)=1$ (see for
instance Theorem 3.38 in \cite{Breiman}.) This completes the proof
of Theorem~\ref{th:recurrence}.

Next we will show that the scaling limit of $X(t)$ is a pure
diffusion.

 For all $n \in {\mathbb N}$ define a sequence of random
continuous increasing functions $T_n(t), 0 \leq t \leq 1$ as
follows. If $t = {m}/{n}$ for some $m$,  $0 \leq m \leq n$ then
\begin{equation}\label{eq:Tnt1}
T_n(t) = \tau_{[nt]}.
\end{equation}
If ${m}/{n} < t < {(m+1)}/{n}$ then $T_n(t)$ is defined by linear
interpolation. That~is,
\begin{equation}\label{eq:Tnt2}
T_n(t) = \tau_m + n (\tau_{m+1} - \tau_m) (t - \frac{m}{n}).
\end{equation}
Let $\{\tau_i\}_{i\geq 0}$ be as in (\ref{eq:taus}). Put
$\alpha=E(\tau_1)$.  We observe that $ X(\tau_{i+1} ) -
X(\tau_i)$, $i\geq 0$,  are i.i.d. random variables with mean zero
and $E(X({\tau_{1}})^2) \leq \lambda
E(\tau_1)+\lambda^2E(\tau_1)^2< \infty $ for all $i\geq 0$. Put
$\beta^2=\var X({\tau_1})$. Note that $\beta>0$. To see this
assume that $A$ is the event that a bond breaks at the first jump
and the broken bond is fixed before the second jump. Then $P(A)=
p\mu/(\lambda+\mu)$, and
\begin{equation*}
P\left(|X_{\tau_1}|>0\right) \geq P\left(|X_{\tau_1}|=1\right)
\geq P\left(A\right)  >0.
\end{equation*}
 Define the sequence of random functions
$\{X_n(\cdot)\}_{n\geq 1}$ by
\[ X_n(t) = \frac{X( n
t)}{\sqrt n },\]  and let $B(t)$, $0 \leq t \leq 1$ denote the
standard one-dimensional Brownian motion. We will prove the
following.
\begin{theorem}\label{th:scaling limit}
\begin{equation*}
X_n(t) \Rightarrow \frac{\beta}{\sqrt{\alpha}}B(t) \quad \mbox{ as
} n \to \infty.
\end{equation*}
\end{theorem}
Let us denote by ${\mathbb P}$ the probability measure on the
underlying probability space $\Omega$. First we prove the
following.

\begin{lemma}\label{lem:taun converging}
\begin{equation*}
\lim_{n \to \infty} \sup_{0 \leq t \leq 1} \left|\frac{T_n(t)}{n}
- \alpha t\right| = 0 \quad{\mathbb P}\mbox{-a.s.}
\end{equation*}
\end{lemma}
{\it Proof } Let $T_n(t)$ be defined by equations (\ref{eq:Tnt1})
and (\ref{eq:Tnt2}). If $t \in(m/n, (m+1)/n)$ then \[\left|T_n(t)
- \tau_{[nt]}\right| \leq (\tau_{m+1}-\tau_m).\] Therefore we have
for all $t \in(m/n, (m+1)/n)$,
\[\frac{1}{n}|T_n(t) - \tau_{[nt]}| \leq
\frac{(\tau_{m+1}-\tau_m)}{n} .\] Using the fact that $\tau_1$ has
finite third moments and the Borel-Cantelli lemma we can easily show
that
\begin{equation}\label{Ttau}
\lim_{n \to \infty} \sup_{0\leq t \leq 1} \frac{1}{n} |T_n(t) -
\tau_{[nt]}| = 0 \quad{\mathbb P}\mbox{-a.s.}
\end{equation}

Let $m \in {\mathbb N}$ be given and consider $t$ of the form
$l/m$, $0 \leq l \leq m$. By the strong law of large numbers we
have for all $l$, $0 \leq l \leq m$
\[
\lim_{n \to \infty} \frac{\tau_{[n \frac{l}{m}]}}{n} = \frac{l}{m}
\alpha \quad{\mathbb P}\mbox{-a.s.}
\]
Using \eqref{Ttau} we have

\[
\lim_{n \to \infty} \frac{T_{n}(\frac{l}{m})}{n} = \frac{l}{m}
\alpha \quad{\mathbb P}\mbox{-a.s.}
\]
Therefore for ${\mathbb P}\mbox{-almost all } \omega \in \Omega$
there exists an $N_m$ such that
\[\left|\frac{T_{n}(\frac{l}{m})}{n}- \frac{l}{m}\alpha\right| <
\frac{\alpha}{m}\quad \mbox{ if $n \geq N_m$.}\]

If $l/{m} < t < {(l+1)}/{m}$, then

\[
\frac{T_n(t)}{n} > \frac{T_{n}(\frac{l}{m})}{n} > \frac{l}{m}
\alpha - \frac{\alpha}{m} = \frac{l+1}{m} \alpha -
\frac{2\alpha}{m}
>\alpha  t - \frac{2 \alpha}{m}.
\]

Similarly
\[
\frac{T_n(t)}{n} < \frac{T_{n}(\frac{l+1}{m})}{n} < \frac{l+1}{m}
\alpha + \frac{\alpha}{m} = \frac{l}{m} \alpha + \frac{2\alpha}{m}
<\alpha  t + \frac{2 \alpha}{m}.
\]

From this it follows that
\[
\lim_{n \to \infty} \sup_{0 \leq t \leq 1} \left|\frac{T_n(t)}{n}
- \alpha t\right| = 0 \quad{\mathbb P}\mbox{-a.s.},
\]
which proves the Lemma.

Now we can prove Theorem~\ref{th:scaling limit}.

{\it Proof of Theorem~\ref{th:scaling limit} } For all $n \in
{\mathbb N} \hbox{ and } 0 \leq t \leq 1$ define $X'_n(t) =
\frac{X(\tau_{[nt]})}{\sqrt n \beta}$. Since $X(\tau_k) =
\sum_{i=1}^k (X(\tau_{i+1} ) - X(\tau_i))$, and $X(\tau_{i+1} ) -
X(\tau_i))$ are i.i.d random variables with finite second moment
it follows from Donsker's invariance principle that
\[
X'_n(t) \Rightarrow B(t) \quad\mbox{ as  $n \to \infty$.}
\]

Let $\{Y(t)\}_{t\geq 0}$ be the process with continuous paths that
is obtained by linearly interpolating between successive jumps of
$\{X(t)\}$. Define $Y_n(t) = \frac{Y(\tau_k)}{\sqrt n \beta}$ if
$t = {k}/{n}$ and if ${k}/{n} < t < {(k+1)}/{n}$, define $Y_n(t)$
by linear interpolation. Since $|\frac{X(\tau_k)}{\sqrt n} -
\frac{Y(\tau_k)}{\sqrt n}| \leq \frac{1}{\sqrt n}$, it follows
that $|X'_n(t)  - Y_n(t)| \leq \frac{1}{\beta \sqrt n}$ for all $
t \in [0,1]$ almost surely.

 This proves that

\[
Y_n(t) \Rightarrow B(t) \quad\mbox{ as  $n \to \infty$.}
\]

For $n\in \mathbb{N}$ and $0\leq t\leq 1$, define
$\lambda_n(t)=T_n(t)/(\alpha n)$. Note that by Lemma~\ref{lem:taun
converging},
\begin{equation*}
\lim_{n\rightarrow \infty} \lambda_n(t)=t\qquad \mbox{uniformly
$\mathbb{P}$-a.s.}
\end{equation*}
Put
 \begin{align*} Y'_n(t)& =\frac{Y(\alpha n
t)}{\sqrt{n} \beta},\\
\intertext{and} Z_n(t) & = Y'_n( \lambda_{n}(t)) =
\frac{Y(T_n(t))}{\sqrt n \beta}.
\end{align*}
Notice that if $t = {k}/{n}$, then $Z_n(t) =
\frac{Y(\tau_k)}{\sqrt n \beta}$.  We will  show that paths of
$Z_n(t)$ and $Y_n(t)$ are uniformly close in $[0,1]$. We will make
use of the fact that all (high enough) moments of $\tau_1$ are
finite and that the jumps occur at Poisson times. The functions
$Y_n(\cdot)$ and $Z_n(\cdot)$ agree at times ${m}/{n}$; while
$Y_n(t)$ linearly interpolates between times ${k}/{n}$ and
${(k+1)}/{n}$, $Z_n(t)$ follows the $Y(t)$ process between times
$\tau_k$ and $\tau_{k+1}$ with rescaled time.  If $t \in
[k/n,(k+1)/n]$, then $\tau_k \leq T_n(t) \leq \tau_{k+1}$, and

\begin{equation*}\sup_{t \in [\frac{k}{n},\frac{k+1}{n}]}\left|\frac{Y(\tau_k)}{\beta
\sqrt n}- Z_n(t)\right| \leq \frac{1}{\sqrt n \beta} \sup_{\tau_k
\leq t\leq \tau_{k+1}}(\left|Y(t)- Y(\tau_k)\right|)
\end{equation*}
and
\begin{equation*}
\sup_{t \in
[\frac{k}{n},\frac{k+1}{n}]}\left|\frac{Y(\tau_k)}{\beta \sqrt n}-
Y_n(t)\right| \leq \frac{1}{\sqrt n \beta} (\left|Y(\tau_{k+1})-
Y(\tau_k)\right|).
\end{equation*}
 Now

\begin{multline*}
{\mathbb P}\left(\sup_{t \in
[\frac{k}{n},\frac{k+1}{n}]}\left|Y_n(t)-
Z_n(t)\right|>\frac{1}{n^{1/4}}\right) \\
\leq {\mathbb P}\left(\sup_{t \in
[\frac{k}{n},\frac{k+1}{n}]}\left|\frac{Y(\tau_k)}{\beta \sqrt n}-
Y_n(t)\right|+ \left|\frac{Y(\tau_k)}{\beta \sqrt n}-
Z_n(t)\right| > \frac{1}{n^{1/4}}\right),
\end{multline*}
which implies
\begin{eqnarray*}
{\mathbb P}\left(\sup_{t \in
[\frac{k}{n},\frac{k+1}{n}]}\left|Y_n(t)-
Z_n(t)\right|>\frac{1}{n^{1/4}}\right) &\leq &{\mathbb
P}\left(\sup_{t\in[\frac{k}{n},\frac{k+1}{n}]}\left|\frac{Y(\tau_k)}{\beta
\sqrt n}- Y_n(t)\right|>\frac{1}{2
n^{1/4}}\right) \\
&+& {\mathbb P}\left(\sup_{t \in
[\frac{k}{n},\frac{k+1}{n}]}\left|\frac{Y(\tau_k)}{\beta \sqrt n}-
Z_n(t)\right| > \frac{1}{2n^{1/4}}\right).
\end{eqnarray*}

Therefore it follows that to show that $Z_n$ and $Y_n$ are
uniformly close it is sufficient to estimate $\sup_{\tau_k \leq
t\leq \tau_{k+1}}(\frac{1}{\sqrt n \beta}|Y(t)- Y(\tau_k)|)$ which
we proceed to do now.

Let

$$M_k = \sup_{\tau_{k-1} \leq t \leq \tau_k} (|Y(t)-
Y(\tau_{k-1})|).$$

We now estimate ${\mathbb P}(M_k > \frac{\beta}{2} n^{1/4})$. Since
the argument for the estimate is the same for all $k$, we estimate
${\mathbb P}(M_1
> \frac{\beta}{2} n^{1/4})$.  Let $C_k$ denote the kth moment of $\tau_1$ and $N(t)$ be defined as in
 the paragraph before the proof of Theorem~\ref{th:recurrence}. We
 have
\begin{eqnarray*}
{\mathbb P}\left(M_1 > \frac{\beta}{2} n^{1/4}\right) & = &
{\mathbb P}\left(M_1
> \frac{\beta}{2} n^{1/4};\tau_1> n^{1/8}\right) + {\mathbb P}\left(M_1
>
\frac{\beta}{2}  n^{1/4}; \tau_1 \leq n^{1/8}\right)\\
& \leq & \frac{C_{17}}{n^{\frac{17}{8}}} + {\mathbb P}\left(N\left(n^{1/8}\right) >\frac{\beta}{2} n^{1/4};\tau_1 \leq n^{1/8}\right)\\
& \leq & \frac{C_{17}}{n^{\frac{17}{8}}} + {\mathbb P}\left(N\left(n^{1/8}\right) > \frac{\beta}{2} n^{1/4}\right)\\
 & \leq & \frac{C_{17}}{n^{\frac{17}{8}}} + A e^{- Bn^{1/8}}
\end{eqnarray*}
for some $A$ and  $B > 0$, where the last inequality follows from
applying the Chebyshev inequality to $\exp(N(n^{1/8}))$. Now
\[
{\mathbb P}\left(\sup_{t \in [0,1]} \left|Z_n(t) - Y_n(t)\right|
> \frac{1}{n^{1/4}}\right) \leq n{\mathbb P}\left(M_1 > \frac{\beta}{2}
n^{1/4}\right).
\]
  Since $\sum_n ( \frac{C_{17}}{n^{\frac{17}{8}}} + A e^{-B
 n^{1/8}}) n < \infty$, it easily follows from Borel-Cantelli lemma
 that
\[
\sup_{0 \leq t \leq 1}\left |Z_n(t) - Y_n(t)\right| \to 0 \quad
{\mathbb P}-a.s.
\]

From this it follows that
\[
Z_n(t) \Rightarrow B(t) \quad\mbox{as  $n \to \infty$.}
\]

 While we have considered the time interval $[0,1]$, we can define all random functions
 discussed so far on the interval $[0,T]$.
 All the above arguments can be easily extended to prove the convergence of
 $Z_n(t)$ to $B(t)$ for $t \in [0,T]$. For the last part of the proof
 we need this convergence to be extended to $[0,T]$, where $T > 1$.
Let us take $T =2$.
We recall that
\[
Y_{n}(\lambda_{n}(t)) = Z_n(t).
\]
By Skorohod's representation theorem, for all $n \in {\mathbb N}$
there exists random variables $\tilde{Y}'_n$, $\tilde{Z}_n$,
$\tilde{B}$, $\tilde{\tau}_n$, and  $\tilde{\lambda}_n$ on some
probability space $(\tilde{\Omega}, \tilde{{\mathbb P}})$, such
that $\tilde{Y}'_n$, $\tilde{Z}_n$, $\tilde{B}$, $\tilde{\tau}_n$,
{ and } $\tilde{\lambda}_n$ have the same distribution
respectively as $Y'_n$, $Z_n$, $B$, $\tau_n$,  and  $\lambda_n$;
$\tilde{Z}_n(t)=\tilde{Y}'_n(\tilde{\lambda}_n(t))$, and both
$\tilde{\lambda}_n(t)\rightarrow t$ and $\tilde{Z}_n(t) \to
\tilde{B}(t)$ uniformly $\tilde{{\mathbb P}}$-almost surly. Thus
$\tilde{B}$ is a $C[0,2]$-valued random variable on
$(\tilde{\Omega}, \tilde{{\mathbb P}})$ with the same distribution
as $B$. Let  $\epsilon>0$ be given and let $\omega \in
\tilde{\Omega}$ be such that
\begin{align}
\lim_{n \to \infty} \tilde{Z}_{n}(t,\omega) & =
\tilde{B}(t,\omega) \quad\mbox{uniformly, } \label{Zcon}\\
\intertext{and}  \lim_{n \to \infty} \tilde{\lambda}_n(t,\omega) &
= t \quad\mbox{uniformly. }\label{lambdacon}
\end{align}
There exists an $N_1= N_1(\omega) \in {\mathbb N}$ such that if $n
\geq N_1$, then $\tilde{\lambda}_n(1, \omega) \in [0,2]$. Take  $n
\geq N_1$. Since $\tilde{\lambda}_n(t, \omega) $ is a continuous
increasing function, for all $t \in [0,1]$, there exists
$e_n(t,\omega) \in [0,2]$ such that
$\tilde{\lambda}_n(e_n(t,\omega), \omega) = t $. Therefore
$\tilde{Y}'_n(t,\omega) = \tilde{Z}_n(e_n(t,\omega),\omega)$. By
\eqref{Zcon} there exists an $N_2 = N_2(\omega) \in {\mathbb N}$,
 $N_2 \geq N_1 $ such that if $n \geq N_2$ then
\begin{equation}\label{eq1}
\sup_{0 \leq t \leq 1} \left|\tilde{Z}_n(e_n(t,\omega),\omega) -
\tilde{B}(e_n(t,\omega),\omega)\right| < \frac{\epsilon}{2}.
\end{equation}
Since $\tilde{B}(t,\omega) \in C([0,2])$ there exists
$\delta(\omega, \epsilon) > 0$ such that if $s$, $s' \in [0,2]$
and $|s-s'|< \delta$ then $|\tilde{B}(s,\omega) -
\tilde{B}(s',\omega)| < {\epsilon}/{2}$.  By \eqref{lambdacon} we
have there exists an $N_3 = N_3(\omega) \in {\mathbb N}$, $N_3
\geq N_2$ such that if $n \geq N_3$, then $\sup_{0 \leq t \leq
1}|t - e_n(t,\omega)| < \delta$. Therefore if $n \geq N_3$,
\begin{equation}\label{eq2}
\left|\tilde{B}(t,\omega) - \tilde{B}(e_n(t,\omega),\omega)\right|
< \frac{\epsilon}{2}.
\end{equation}

From equations \eqref{eq1} and \eqref{eq2} we have that if $n \geq
N_3$, for all $t \in [0,1]$,
\begin{align*}
 \left|\tilde{Y}'_n(t,\omega) - \tilde{B}(t,\omega)\right| &{} =
\left|\tilde{Z}_n(e_n(t,\omega),\omega) - \tilde{B}(t,\omega)\right|&\\
{}&{}\leq \left|\tilde{Z}_n(e_n(t,\omega),\omega) -
\tilde{B}(e_n(t,\omega),\omega)\right|& \\
{}&{}\quad+\left|\tilde{B}(e_n(t,\omega),\omega)-\tilde{B}(t,\omega)
\right|<  \epsilon.
\end{align*}

This proves that
\[
\tilde{Y}'_n(t) \to \tilde{B}(t)\quad \mbox{uniformly in $[0,1]$ }
\tilde{{\mathbb P}}\mbox{-a.s.},
\]
which implies that
\[
Y'_n(t) \Rightarrow B(t) \quad\mbox{as  $n \to \infty$.}
\]

Since $|X(t)-Y(t)|\leq 1$ for all $t\geq 0$,  we have that
\begin{equation*}
\left|\frac{X(\alpha nt)}{\sqrt{n}}-\frac{Y(\alpha n
t)}{\sqrt{n}}\right|\rightarrow 0 \qquad\mbox{uniformly as
$n\rightarrow \infty$},
\end{equation*}
and thus
\[ \frac{X(\alpha n t)}{\beta \sqrt{n}} \Rightarrow B(t),\]
proving the theorem.

The following proposition establishes a lower bound for
$\beta/\sqrt{\alpha}$ in  some cases.
\begin{proposition}\label{pr:variance identity}
 If $p=1$ and $\mu$ is large enough, then
\begin{equation}\label{eq:variance identity}
E\left(X(\tau_1)^2\right)>\lambda E\left(\tau_1\right).
\end{equation}
\end{proposition}
{\it Proof } Since $p=1$, the first jump happens at $\sigma_1$.
Thus $P(X(\sigma_1)=1)=P(X(\sigma_1)=-1)=1/2$. Let $\zeta$ be the
first time after $\sigma_1$ that either the particle jumps or the
broken bond is fixed. Therefore
\begin{align*}
E\left(\zeta-\sigma_1\right)& = \frac{1}{\lambda+\mu}\\
\intertext{and}
 E\left(\zeta\right)&=\frac{1}{\lambda}+ \frac{1}{\lambda+\mu}.
\end{align*}

Note that the probability that the broken bond is fixed before the
particle jumps is $\mu/(\lambda+\mu)$. Therefore,
\begin{align*}
P\left(X(\zeta)=1\mid X(\sigma_1)=1\right) & =
\frac{\mu}{\lambda+\mu},\\
\intertext{and} P\left(X(\zeta)=2 \mid X(\sigma_1)=1\right) & =
\frac{\lambda}{\lambda+\mu}.
\end{align*}
Hence we have
\begin{equation*}
E\left(X(\zeta)^2\right)=1^2\cdot\frac{\mu}{\lambda+\mu}+2^2\cdot\frac{\lambda}{\lambda+\mu}=
\frac{\mu+4\lambda}{\lambda+\mu}
\end{equation*}
and
\begin{equation}\label{eq:variance of zeta}
E\left(X(\zeta)^2\right)-\lambda E\left(\zeta\right) =
\frac{2\lambda}{\lambda+\mu}>0.
\end{equation}

By definition of $\zeta$ and $\tau_1$,
$P(\tau_1=\zeta)=\mu/(\lambda+\mu)$.  On $\{\tau_1>\zeta\}$, there
are two broken bonds at time $\zeta$. Given that
$\{\tau_1>\zeta\}$, the conditional probability that both these
bonds get fixed before there is a third jump is
\begin{equation*}
\frac{2\mu}{\lambda+2\mu}\cdot\frac{\mu}{\lambda+\mu}=1-O\left(\frac{1}{\mu}\right).
\end{equation*}
Therefore, the chance that there is a third jump is of the order
$1/{\mu^2}$. Hence, as $\mu\rightarrow \infty$, the quantity
$1-P(X(\tau_1)=X(\zeta))$ is at most of the order $1/{\mu^2}$, and
\begin{equation}\label{eq:compare variances}
E\left(X(\tau_1)^2\right)+O\left(\frac{1}{\mu^2}\right)\geq
E\left(X(\zeta)^2\right) \mbox{ as $\mu\rightarrow \infty$.}
\end{equation}

Let $\theta$ be the time of the first jump after $\zeta$. Note
that $\{\tau_1>\theta\}\subset\{\tau_1>\zeta\}$ and therefore,
$P(\tau_1>\theta)$ is at most of order $1/{\mu^2}$. Furthermore,
$E(\theta-\zeta\mid \tau_1>\theta)=1/{\lambda}$ and so
$E(\theta-\zeta)$ is at most of order $1/{\mu^2}$. Hence,
$E(\min(\tau_1-\zeta,\theta-\zeta))$ is at most of order
$1/{\mu^2}$. Also, using a first step analysis, we can show that
for all $\mu\geq 1$,  $E(\tau_1-\theta\mid \tau_1>\theta)\leq
4(e^{\lambda}-1)/\lambda$, and therefore $E(\tau_1-\theta)$ is
 of order $1/{\mu^2}$. Since
\begin{equation*}
\tau_1-\zeta\leq \min(\tau_1-\zeta,\theta-\zeta)+ (\tau_1-\theta),
\end{equation*}
$E(\tau_1-\zeta)$ is also of  order at most $1/{\mu^2}$. Along
with (\ref{eq:variance of zeta}) and (\ref{eq:compare variances}),
this establishes the Proposition.

Without loss of generality, for the rest of this discussion,  we
will assume that $\lambda=1$. Note that if $U(t)$ is a continuous
time nearest neighbor random walk with jump rate 1, then
\begin{equation*}
\frac{U(nt)}{\sqrt{n}}\Rightarrow B(t).
\end{equation*}
On the other hand, Theorem~\ref{th:scaling limit} and
Proposition~\ref{pr:variance identity} imply that under the
conditions of Proposition~\ref{pr:variance identity}, we have
\begin{equation*}
\frac{X(nt)}{\sqrt{n}} \Rightarrow
\frac{\beta}{\sqrt{\alpha}}B(t),
\end{equation*}
with $\beta/\sqrt{\alpha}>1$. Thus, under the conditions of
Proposition~\ref{pr:variance identity}, the effect of bond
breaking is to make the process move faster in the scaling limit.

\bigskip

{\bf Remark:}  We can prove the scaling limit for the
corresponding walk in $d$ dimensions by obtaining the bound on the
recurrence time $\tau$ in the same way using an $M/M/\infty$ queue
and in the proof of Theorem 1 using the multidimensional version
of Donsker's invariance principle. The rest of the argument
proceeds in the same way.

{\bf Acknowledgements:} Majid Hosseini would like to thank Burgess
Davis for helpful discussions. Ravishankar would like to thank
Raghu Varadhan for interesting discussions and suggestions.


\end{document}